\documentclass[11pt]{article}
\usepackage{anysize}
\marginsize{2.0cm}{2.0cm}{2.0 cm}{2.0 cm}
\usepackage{setspace}

\usepackage{latexsym,amsfonts,amsmath}
\usepackage{graphicx}
\usepackage{color}
\usepackage{epstopdf}
\usepackage[titletoc]{appendix}

\def\RR{\mathbb{R}}

\def\II{\mathbb{I}}

\def\liminf{\mathop{\underline{\lim}}}

\newcommand{\defi}{:=}

\newtheorem{condition}{Condition}[section]{\bfseries}{\itshape}

{\bfseries}{\itshape}

{\bfseries}{\itshape}

\newtheorem{theorem}{Theorem}[section]{\bfseries}{\itshape}

\newtheorem{corollary}{Corollary}[section]{\bfseries}{\itshape}

\newtheorem{proposition}{Proposition}[section]{\bfseries}{\itshape}

{\bfseries}{\itshape}

\newtheorem{lemma}{Lemma}[section]{\bfseries}{\itshape}

\newtheorem{remark}{Remark}[section]{\bfseries}{\itshape}

\newtheorem{definition}{Definition}[section]{\bfseries}{\itshape}

{\bfseries}{\itshape}

\begin{document}
\title{Linear programming approach to optimal impulse control problems with functional constraints\thanks{Declarations of interest: none.}}
\date{}
\author{
Alexey Piunovskiy\thanks{Department of Mathematical Sciences, University of Liverpool, L69 7ZL, UK. E-mail: piunov@liv.ac.uk}~ and
Yi Zhang\thanks{Corresponding author.
Department of Mathematical Sciences, University of Liverpool, L69 7ZL, UK. E-mail: yi.zhang@liv.ac.uk}
}

\maketitle

\par\noindent\textbf{Abstract}.
This paper considers an optimal impulse control problem of dynamical systems generated by a flow. The performance criteria are total costs over the infinite time horizon. Apart from the main performance to be minimized, there are multiple constraints on performance functionals of a similar type. Under a natural set of compactness-continuity conditions on the system primitives, we establish a linear programming approach, and prove the existence of a stationary optimal control strategy out of a more general class of randomized strategies. This is done by making use of the tools from Markov decision processes.

\bigskip


\section{Introduction}

The literature of impulse control problems is quite vast in terms of both theoretical developments and applications. One popular method of investigations is dynamic programming, based on solving the quasi-variational inequality, see e.g., \cite{Barles:1985,Bensoussan:1975,DufourHoriguchiPiunovskiy:2016,Dufour:2015,b12,Miller:2020,PiunovskiySasha:2018}, where the system under control is deterministic in \cite{Barles:1985,PiunovskiySasha:2018} and stochastic in the others. Various meaningful problems arising from Internet congestion control, finance, epidemiology, inventory management, and insurance were solved using this method in \cite{b9,Korn:1997,PiunovskiySasha:2018,Presman:2006,Wei:2010}. Below we only focus on the literature on deterministic impulse control problems. Another method is the maximum principle. This was established for deterministic impulse control problems in \cite{Blaquiere:1985,Chahim:2012,b8,Miller:1992,Rempala:1988}, where in \cite{b8,Miller:1992}, the original impulse control problems were first transformed to equivalent problems without impulses in discrete-time and in continuous-time, respectively, and then the versions of maximum principle were formulated for the transformed problems. The control model in \cite{Miller:1992} was further generalized in \cite{Miller:2009}, see also the monograph \cite{Miller:2003}. In relation to the maximum principle, certain impulse control problems with a fixed number of impulses can be viewed as mathematical programming problems with the collection of admissible impulse moments and impulses being the control variables. This method was demonstrated in \cite{HouWong:2011,Liu:1998}, where the gradient of the performance functionals with respect to the control variables was calculated.

A different method from the aforementioned ones is the linear programming approach. This method is rather powerful in handling problems with functional constraints. The formulation of the linear program is often based on the observation that the performance can be expressed as a linear function of suitably defined occupation measures, which satisfy a certain linear equation. Then one may formulate the linear program on the space of all measures satisfying that linear equation. The linear programming approach is established if one can show that the linear program and the original impulse control problem are equivalent. In general, this equivalence is a delicate matter, and its validity depends on the nature of the problems. For deterministic optimal control problem without impulses, the linear programming method was investigated and established in e.g., \cite{Gaitsgory:2009,HernandezHernandez:1999,Lasserre:2008}. Similar approach to gradually controlled stochastic systems was justified in e.g., \cite{Bhatt:1996,Buckdahn:2011,Kurtz:1998}.

In the literature on deterministic impulse control problems, to the best of our knowledge, only few works were devoted to the linear programming approach: see \cite{b13,b13p,b13pp}. In \cite{b13} (resp.,\cite{b13p,b13pp}),  for a gradual-impulse control problem over a finite horizon, the authors introduced suitable occupation measures (resp.,made use of tools from the theory of partial differential equations) to formulate linear programs on the space of finite measures. However, the equivalence between the linear program and the original problem was left open in \cite{b13}, and was discussed in \cite{b13p,b13pp}. In the aforementioned works, the flow in the Euclidean space came from an ordinary differential equation, whereas in the present paper the flow is arbitrary enough and lives in a Borel space. Moreover, we allow a simultaneous sequence  of impulses which was not the case in \cite{b13,b13p,b13pp}.

The contributions of the present paper are twofold. First, we establish fully the linear programming approach to a deterministic impulse control problem over an infinite time horizon with total cost criteria, where apart from the main performance to be maximized, there are multiple constraints on other performance functionals of a similar type. We achieve this by transofrming the impulse control problem to an equivalent Markov decision process (MDP) with total cost, and exploit the relevant results obtained for MDPs in \cite{Dufour:2012}. For this reason, the occupation measures introduced here are different from  those introduced in \cite{b13,b13p,b13pp}. Also, since we allow accumulation of impulses, the linear program is on the space of infinite-valued measures.  In Remark \ref{rema} we introduce `aggregated occupation measures', which have  more in common with the occupation measures in \cite{b13,b13p,b13pp}. The corresponding linear programs  in terms of such measures can be found in the follow-up paper \cite{PiunovskiyZhang:2020II}, where we compare our approach with \cite{b13,b13p,b13pp}.

The second contribution of this paper is that we prove the solvability of the linear program and the original impulse control problem, under a natural set of compactness-continuity conditions on the system primitives. For problems with functional constraints, it is natural to consider randomized control strategies, see \cite{Piunovskiy:1997}. We show that there exists a stationary  control strategy, which is optimal in that class. It turns out that the imposed compactness-continuity conditions on the system primitives of the impulse control problem, while being natural, do not imply the equivalent MDP to be a semicontinuous model. To get over this difficulty, we introduce an auxiliary MDP with an extended state space, and, by introducing a suitable metric on it, we show that the auxiliary MDP model is semicontinuous. From this, we eventually retrieve the desired solvability result for the original impulse control problem.

The rest of this paper is organized as follows. In Section \ref{sec2} we state the concerned deterministic impulse control problem under consideration. In Section \ref{15JuneSec02} we establish the linear programming approach for the impulse control problem by reformulating it as a constrained MDP. In Section \ref{sec3}, we obtain the main solvability result for the impulse control problem. We finish this paper with a conclusion in Section \ref{15JuneSectionConclusion}. To improve the readability, some technical results are collected in the appendix.

\section{Problem Statement}\label{sec2}

In what follows, $\II\{\cdot\}$ denotes the indicator function.

We will deal with an impulse control model $\{\textbf{X},\textbf{A},\phi, l, \{(C_j^g,C_j^I)\}_{j=0}^J\}$ defined through the following system primitives:
\begin{itemize}
\item $\bf X$ is the state space, which is a Borel space, i.e., a topological space that is homeomorphic to a Borel subset of a complete separable metric (i.e., Polish) space. We endow $\textbf{X}$ with its Borel $\sigma$-algebra ${\cal B}(\textbf{X}).$
\item $\phi:{\bf X}\times [0,\infty)\to{\bf X}$ is the measurable flow possessing the semigroup property $\phi(x,t+s)=\phi(\phi(x,s),t)$ for all $x\in{\bf X}$ and $(t,s)\in  [0,\infty)\times  [0,\infty)$; $\phi(x,0)=x$ for all $x\in{\bf X}$. Between the consecutive impulses, the state evolves according to the flow.
\item $\bf A$ is the action space, which is a Borel space and endowed with its Borel $\sigma$-algebra ${\cal B}(\textbf{A}).$
\item $l:{\bf X}\times{\bf A}\to{\bf X}$ is the mapping describing the new state after the corresponding impulse is applied at the current state.
\item For each $j=0,1,\dots,J,$ where $J$ is a fixed natural number, $C^g_j:{\bf X}\to [0,\infty)$ is the gradual (or say running) cost rate.
\item For each $j=0,1,\dots,J,$ $C^I_j:{\bf X}\times{\bf A}\to [0,\infty)$ is the cost function associated with the impulses applied at the corresponding states.
\end{itemize}
All the mappings $\phi,l,\{C^g_j\}_{j=0}^J$ and $\{C^I_j\}_{j=0}^J$ are assumed to be (Borel) measurable in their arguments.

Let us describe verbally the system dynamics under control as follows.
If the current state is $x_0\in\textbf{X},$ the decision maker chooses a pair $(\theta_1,a_1)\in [0,\infty]\times \textbf{A}$, where $\theta_1$ represents the time to wait until the next impulse is applied, and $a_1$ is the next impulse. Over $[0,\theta_1)$, the state process evolves according to $\phi(x_0,t)$, and at $\theta_1,$ the state is $l(\phi(x_0,\theta_1),a_1)$. Then the decision maker chooses $(\theta_2,a_2)$, and so on. Note that a simultaneous sequence of impulses is not forbidden: $\theta_i=\theta_{i+1}=\ldots\theta_{k}=0$, and at the corresponding time moment the process has $2+k-i$ different values.

If $\theta_1,\theta_2,\dots,\theta_{i-1}<\infty$,  at the current state $x_{i-1},$ the total cost accumulated over the next interval of length $\theta_i$ equals
\begin{equation}\label{e1pprim}
\int_{[0,\theta_i) } C_j^g(\phi(x_{i-1},u))du+\II\{\theta_i<\infty\} C^I_j(\phi(x_{i-1},\theta_i),a_i).
\end{equation}
In the previous sum, the last summand is absent when $\theta_i=\infty$ because if $\theta_i=\infty$ is selected, it means waiting indefinitely until the next impulse moment, i.e., no further impulse will be applied. When $\theta_i=\infty,$ we will artificially put the next state $x_i=\Delta$ for some isolated point $\Delta\notin \textbf{X}$. We put ${\bf X}_\Delta:={\bf X}\cup\{\Delta\}$. On the other hand, the selection of $\theta_i=0$ means an instantaneous application of the impulse $a_i$ at the current state $x_{i-1}$, leading to the next state $x_i=l(x_{i-1},a_i).$

Mathematically, the system dynamics (also called trajectories henceforth) can be represented as sequences in one of the following two forms:
\begin{eqnarray}
&&x_0\to (\theta_1,a_1)\to x_1\to (\theta_2,a_2)\to \ldots;~~~~ \theta_i<\infty \mbox{ for all } i=1,2,\dots, \nonumber\\
\mbox{or}&&\label{e1}\\
&&x_0\to (\theta_1,a_1)\to\ldots\to x_n\to (\infty,a_{n+1})\to \Delta\to (\theta_{n+2},a_{n+2}) \to \Delta \to \ldots,\nonumber
\end{eqnarray}
where $x_0\in{\bf X}$ is the initial state of the controlled process and $\theta_i<\infty$ for all $i=1,2,\ldots, n$.

Let the space of all the trajectories (\ref{e1}) be denoted by
\begin{eqnarray*}
\Omega&=&\bigcup_{n=1}^\infty[{\bf X}\times(([0,\infty)\times{\bf A})\times{\bf X})^n
\times(\{\infty\}\times{\bf A})\times\{\Delta\}\times (([0,\infty]\times{\bf A})\times\{\Delta\})^\infty]\\
&& \cup [{\bf X}\times(([0,\infty)\times{\bf A})\times{\bf X})^\infty].
\end{eqnarray*}
We fix the natural Borel $\sigma$-algebra $\cal F$ on $\Omega$.
Finite sequences in the form
\begin{eqnarray*}
h_i=(x_0, (\theta_1,a_1), x_1, (\theta_2,a_2),\ldots,x_i)
\end{eqnarray*}
will be called $i$-histories; $i=0,1,2,\ldots$. The space of all such $i$-histories will be denoted as ${\bf H}_i$. We endow it with the $\sigma$-algebra ${\cal F}_i:={\cal B}({\bf H}_i)$, which is the restriction of ${\cal F}$ to ${\bf H}_i$.

In line with the previous verbal description, for $i=1,2,\dots,$ at the current state $x_{i-1}\in{\bf X}$, if the pair $(\theta_i,a_i)\in[0,\infty]\times{\bf A}$ is the selected control, then after $\theta_i$ time units, the impulse $a_i$ will be applied, leading instantaneously to the new state
\begin{equation}\label{e1p}
x_i=\left\{\begin{array}{ll}
l(\phi(x_{i-1},\theta_i),a_i), & \mbox{ if } \theta_i<\infty; \\
\Delta, & \mbox{ if } \theta_i=\infty.
\end{array}\right.
\end{equation}
Once $\Delta$ appeared for the first time, it will appear indefinitely in the system dynamics, and one may regard it as an absorbing state.

Next, let us describe the selection of the pairs $(\theta_i,a_i)$. If $\theta_1,\theta_2,\ldots,\theta_{i-1}<\infty$, at the current state $x_{i-1}$, the decision maker has in hand the information about the $i-1$-history, that is, the sequence
$
h_{i-1}=(x_0, (\theta_1,a_1), x_1,\ldots,  (\theta_{i-1},a_{i-1}),x_{i-1}).
$
The selection of the next pair $(\theta_i,a_i)\in[0,\infty]\times \textbf{A}$ is based on this information, and we further allow the selection of the pair $(\theta_i,a_i)$ to be randomized. For this reason, even though the original uncontrolled process evolves deterministically according to the flow $\phi$, a suitable probability space will be constructed.

Control strategies (or simply say strategies) specify the selection of the pair $(\theta_i,a_i)$ given the $i-1$-history $h_{i-1}$, and are defined as follows.

\begin{definition}\label{d1}
\begin{itemize}
\item[(a)]
A (control) strategy $\pi=\{\pi_i\}_{i=1}^\infty$ is a sequence of stochastic kernels $\pi_i$ on $[0,\infty]\times\textbf{A}$ given $\textbf{H}_{i-1}$.
\item[(b)]
A Markov strategy is defined by the sequence of stochastic kernels in the form $\{\pi_i(d\theta\times da|x_{i-1})\}_{i=1}^\infty$.
\item[(c)]
A strategy is called stationary and denoted as $\widetilde{\pi}$, if there is a stochastic kernel $\widetilde{\pi}$ on $[0,\infty]\times\textbf{A}$ given $\textbf{X}_\Delta$  such that $\pi_i(d\theta\times da|h_{i-1})=\widetilde{\pi}(d\theta\times da|x_{i-1})$ for all $i=1,2,\ldots$.
\item[(d)]
A strategy ${\pi}=\{\pi_i\}_{i=1}^\infty$ is called deterministic stationary and denoted as $f$, if for all $i=1,2,\ldots$, $\pi_i(d\theta\times da|h_{i-1})=\delta_{f(x_{i-1})}(d\theta\times da)$, where $f:\textbf{X}_\Delta\to[0,\infty]\times  \textbf{A}$ is a measurable mapping.
\end{itemize}
\end{definition}
In terms of interpretation, under a strategy $\pi=\{\pi_i\}_{i=1}^\infty$, $\pi_i(d\theta\times da|h_{i-1})$ is the (regular) conditional distribution of $(\Theta_{i},A_i)$\footnote{As was mentioned above, the selection of $(\theta_i,a_i)$ can be randomized, so that any selected $(\theta_i,a_i)$ can be regarded as the realization of the pair of some random variables $(\Theta_{i},A_i)$.} given the $i-1$-history $H_{i-1}=h_{i-1}$. This is in line with the following construction. Here and below, the capital letters $X_i,T_i,\Theta_i, A_i$ and $H_i$ denote the corresponding functions of $\omega\in\Omega$, i.e., random elements.

For a given initial distribution $\nu$ on $\textbf{X}$ and a strategy $\pi$, by the Ionescu-Tulcea Theorem, see e.g., \cite[Prop.7.28]{Bertsekas:1978}, there is a unique probability measure $P^\pi_{\nu}$ on $(\Omega,{\cal F})$ satisfying the following conditions:
\begin{eqnarray*}
P^\pi_{\nu}(X_0\in\Gamma_X)=\nu(\Gamma_X)~\forall~ \Gamma_X\in{\cal B}({\bf X}_\Delta);
\end{eqnarray*}
and for all $i=1,2,\dots$, $\Gamma\in{\cal B}([0,\infty]\times{\bf A})$, $\Gamma_X\in{\cal B}({\bf X}_\Delta)$,
\begin{eqnarray}
P^\pi_{\nu}((\Theta_i,A_i)\in\Gamma|H_{i-1})&=&\pi_i(\Gamma |H_{i-1});\label{e33}\\
P^\pi_{\nu}(X_i\in\Gamma_X|H_{i-1},(\Theta_i,A_i))
&=& \left\{\begin{array}{ll}
\delta_{l(\phi(X_{i-1},\Theta_i),A_i)}(\Gamma_X), & \mbox{ if } X_{i-1}\in{\bf X},~\Theta_i<\infty; \\
\delta_\Delta(\Gamma_X) & \mbox{ otherwise. } \end{array} \right. \nonumber
\end{eqnarray}
When the initial distribution $\nu$ is a Dirac measure concentrated on a singleton, say $\{x_0\}$, we write $P^\pi_\nu$ as $P^\pi_{x_0}$. The mathematical expectation with respect to $P^\pi_\nu$ and $P^\pi_{x_0}$ is denoted as $E^\pi_\nu$ and $E^\pi_{x_0}$, respectively.

Let us introduce the performance measure of a strategy $\pi$ at an initial distribution $\nu$:
\begin{eqnarray*}
{\cal V}_j(\nu, \pi) &:=&  E^\pi_{\nu}\left[\sum_{i=1}^\infty \II\{X_{i-1}\ne\Delta\} \left\{ \int_{[0,\Theta_{i})}  C^g_j(\phi(X_{i-1},u)) du\right.\right.+ \II\{\Theta_i<\infty\} \left.\left.\vphantom{\sum_{i=1}^\infty} C^I_j(\phi(X_{i-1},\Theta_i),A_i)\right\}\right]
\end{eqnarray*}
for each $j=0,1,\dots,J.$ Again, when $\nu=\delta_{x_0}$, we write ${\cal V}_j(x_0, \pi)$ for ${\cal V}_j(\nu, \pi)$. Note that we do not exclude the possibility of $\sum_{i=1}^{\infty}\Theta_i<\infty$, but we will only consider the total cost accumulated over $[0,\sum_{i=1}^\infty \Theta_i).$ This is consistent with the definition of ${\cal V}_j(\nu, \pi)$.

The constrained optimal control problem under study is the following one:
\begin{eqnarray}\label{PZZeqn02}
\mbox{Minimize with respect to } \pi &&{\cal V}_0(x_0, \pi)  \\
\mbox{subject to }&& {\cal V}_j(x_0,\pi)\le d_j,~j=1,2,\dots,J.\nonumber
\end{eqnarray}
Here and below, we take $x_0\in\textbf{X}$ as a fixed initial point, and $\{d_j\}_{j=1}^J$ as fixed constraint constants.

\begin{definition}\label{d2}
A strategy $\pi$ is called feasible if it satisfies all the constraint inequalities in problem (\ref{PZZeqn02}). A feasible strategy $\pi^\ast$ is called optimal if ${\cal V}_0(x_0,\pi^*)\le  {\cal V}_0(x_0,\pi)$ for all feasible strategies $\pi$.
\end{definition}

The primary goal of this paper is to show, under a natural set of conditions, the existence of a stationary optimal strategy for the impulse control problem (\ref{PZZeqn02}), and to establish the linear programming approach that can be used to obtain it. It is thus outside our concern if either problem (\ref{PZZeqn02}) has no feasible strategy or all the feasible strategies $\pi$ are with an infinite value, i.e., ${\cal V}_0(x_0,\pi)=\infty$. In the former case, problem (\ref{PZZeqn02}) is not solvable; and in the latter case, any feasible strategy is optimal. Therefore, in the forthcoming discussions, we shall assume that the following condition is in force, regarding the consistency of problem (\ref{PZZeqn02}).
\begin{condition}\label{ConstrainedPPZcondition05}
There exists some feasible strategy $\pi$ such that ${\cal V}_0(x_0,\pi)<\infty.$
\end{condition}

The solvability of problem (\ref{PZZeqn02}) cannot be guaranteed without imposing further requirements on the system primitives. We shall impose the following set of compactness-continuity conditions, under which, we will actually show the existence of a stationary optimal strategy, and establish the linear programming approach for obtaining it.
\begin{condition}\label{ConstrainedPPZcondition01}
\begin{itemize}
\item[(a)] The space $\textbf{A}$ is compact, and $\infty$ is the one-point compactification of $[0,\infty)$.
\item[(b)] The mapping $(x,a)\in \textbf{X}\times \textbf{A}\rightarrow l(x,a)$ is continuous.
\item[(c)] The mapping $(x,\theta)\in \textbf{X}\times [0,\infty)\rightarrow \phi(x,\theta)$ is continuous.
\item[(d)] For each $j=0,1,\dots,J,$ the function $(x,a)\in \textbf{X}\times \textbf{A}\rightarrow C_j^I(x, a)$ is lower semicontinuous.
\item[(e)] For each $j=0,1,\dots,J,$ the function $x\in \textbf{X}\rightarrow C_j^g(x)$ is lower semicontinuous.
\end{itemize}
\end{condition}

\section{Reformulation as an MDP and as a linear program}\label{15JuneSec02}
In this section, we reformulate the impulse control problem (\ref{PZZeqn02}) as an MDP with total cost criteria. By exploiting some facts about this class of MDPs, we will establish the relation between the impulse control problem (\ref{PZZeqn02}) and a linear program. While the MDP formulation itself is natural, the resulting MDP is not always convenient to deal with: under Condition \ref{ConstrainedPPZcondition01}, it is not a semicontinuous model. We shall elaborate on and deal with it in the next section (see the proof of Theorem \ref{Theorem01}    therein), where the main issue of solvability will be addressed.

\subsection{The MDP formulation}
The formulation of the impulse control model $\{\textbf{X},\textbf{A},\phi, l, \{(C_j^g,C_j^I)\}_{j=0}^J\}$ as an MDP is done as follows. The state space of the MDP is
$
{\bf X}_\Delta=\textbf{X}\cup \{\Delta\},
$
where $\Delta\notin \textbf{X}$ is an isolated point. The action space is
$
\textbf{B}:=[0,\infty]\times{\bf A},
$
which is endowed with the product topology and the corresponding Borel $\sigma$-algebra. The transition kernel is defined by
\begin{eqnarray*}
Q(dy|x,(\theta,a)):=\left\{\begin{array}{ll}
\delta_{l(\phi(x,\theta),a)}(dy), & \mbox{ if } x\ne\Delta,~\theta\ne \infty;\\ \delta_\Delta(dy) & \mbox{ otherwise}, \end{array}\right..
\end{eqnarray*}
The cost functions are given by
\begin{eqnarray*}
\bar{C}_j(x,(\theta,a))= \II\{x\ne\Delta\} \left\{ \int_{[0,\theta]}  C^g_j(\phi(x,u)) du
+ \II\{\theta<\infty\} C^I_j(\phi(x,\theta),a)\right\},~j=0,1,\dots,J,
\end{eqnarray*}
and the constraint constants are $d_j\in[0,\infty),~j=1,2,\dots,J.$ Here $J$ is the number of constraints. We may summarize this MDP model as $\{\textbf{X}_\Delta,\textbf{B},Q,\bar{C}_0,\{\bar{C}_j,d_j\}_{j=1}^J\}$.  We shall often refer to this MDP as the ``original'' model, since we will introduce an ``auxiliary'' MDP model in the proof of Theorem \ref{Theorem01} below.

Now the impulse control problem (\ref{PZZeqn02}) can be viewed as a (constrained) MDP with total cost criteria.

\begin{definition}\label{15JuneDef02}
For each strategy $\pi$, its occupation measure $\mu^\pi$ in the original MDP model is defined by
\begin{eqnarray}\label{enum1}
\mu^\pi(\Gamma_1\times\Gamma_2):=E_{x_0}^\pi\left[\sum_{n=0}^\infty \II\{(X_n,B_{n+1})\in\Gamma_1\times\Gamma_2\}\right]~\forall~\Gamma_1\in{\cal B}(\textbf{X}_\Delta),\Gamma_2\in{\cal B}(\textbf{B}).
\end{eqnarray}
\end{definition}
Note that the measure $\mu^\pi$ is $[0,\infty]$-valued.

The following fact concerning the occupation measure of a stationary strategy is frequently used below.
\begin{remark}\label{15JuneRemark01}
If $\widetilde{\pi}$ is a stationary strategy, then $\mu^{\widetilde{\pi}}(dx\times db)=\widetilde{\pi}(db|x)\mu^{\widetilde{\pi}}(dx\times \textbf{B})$. This follows from $\mu^\pi(\Gamma_1\times\Gamma_2):=\sum_{n=0}^\infty  E_{x_0}^\pi\left[\II\{(X_n\in \Gamma_1\} E_{x_0}^\pi[  \II\{B_{n+1}\in\Gamma_2\}|X_n]\right]$ in the notation of Definition \ref{15JuneDef02}. (For an arbitrary strategy $\pi$, the decomposition of its occupation measure into the product of its marginal and a stochastic kernel on $\textbf{B}$ given $\textbf{X}_\Delta$ is not trivial because the occupation measure in an MDP with total cost is not finite-valued. See Proposition \ref{15JuneProposition01} below.)
\end{remark}

Let us consider the following linear program for the MDP $\{\textbf{X}_\Delta,\textbf{B},Q,\bar{C}_0,\{\bar{C}_j,d_j\}_{j=1}^J\}$:
\begin{eqnarray}\label{SashaLp01}
\mbox{Minimize over all measures $\mu$ on ${\bf X}_\Delta\times{\bf B}$} &:&\int_{\textbf{X}_\Delta \times\textbf{B}}\bar{C}_0(x,b)\mu(dx\times db) \\
 \mbox{subject to}&:& \mu(dx\times \textbf{B})=\delta_{x_0}(dx)+\int_{\textbf{X}_\Delta\times \textbf{B}}Q(dx|y,b)\mu(dy\times db);\nonumber\\
&&\int_{\textbf{X}_\Delta \times\textbf{B}}\bar{C}_j(x,b)\mu(dx\times db)\le d_j,~j=1,2,\dots,J.\nonumber
\end{eqnarray}

\begin{definition}
\begin{itemize}
\item[(a)]
A feasible measure $\mu$ in linear program (\ref{SashaLp01}) is called to be with a finite value if it satisfies
$
\int_{\textbf{X}_\Delta\times\textbf{B}}\bar{C}_0(x,b)\mu(dx\times db)<\infty.
$
\item[(b)] In the linear program (\ref{SashaLp01}), a measure $\mu_1$ is said to outperform another measure $\mu_2$  if
\begin{eqnarray*}
\int_{\textbf{X}_\Delta\times\textbf{B}}\bar{C}_j(x,b)\mu_1(dx\times db)\le \int_{\textbf{X}_\Delta\times\textbf{B}}\bar{C}_j(x,b)\mu_2(dx\times db)
\end{eqnarray*}
for each $j=0,1,\dots,J.$
\end{itemize}
\end{definition}

Since we may write ${\cal V}_j(x_0,\pi)=\int_{{\bf X}_\Delta\times{\bf B}} \bar{C}_j(x,b)\mu^\pi(dx\times db)$, the occupation measure of any feasible strategy $\pi$ for problem (\ref{PZZeqn02}) satisfies the constraint inequalities in the linear program (\ref{SashaLp01}). Moreover, by \cite[Lem. 9.4.3]{HernandezLerma:1999}, it also satisfies the constraint equality in (\ref{SashaLp01}). Thus, for each feasible strategy $\pi$, its occupation measure $\mu^\pi=\mu$ is feasible for the linear program (\ref{SashaLp01}). Moreover, under Condition \ref{ConstrainedPPZcondition05}, there is some feasible measure with a finite value for the linear program (\ref{SashaLp01}). We may draw from these observations a first relation between the impulse control problem (\ref{PZZeqn02}) and the linear program (\ref{SashaLp01}):
\begin{eqnarray}\label{15JuneEqn05}
\mbox{optimal value of problem (\ref{PZZeqn02})}\ge~ \mbox{optimal value of the linear program (\ref{SashaLp01})}.
\end{eqnarray}

In the opposite direction, under Conditions \ref{ConstrainedPPZcondition05} and \ref{ConstrainedPPZcondition01}, any feasible measure with a finite value for the linear program (\ref{SashaLp01}) is outperformed by the occupation measure of some feasible strategy for problem (\ref{PZZeqn02}). The exact and more informative formulation of this fact, which is given in Proposition \ref{15JuneProposition01} below, uses the following notation. Let
\begin{eqnarray}\label{e5p}
V:=\left\{x\in\textbf{X}_\Delta: \inf_{\pi}E_x^\pi\left[ \sum_{n=0}^\infty \sum_{j=0}^J\bar{C}_j(X_n,B_{n+1}) \right]>0\right\},
\end{eqnarray}
where $\{B_n\}_{n=1}^\infty$ is the action process in the MDP $\{\textbf{X}_\Delta,\textbf{B},Q,\bar{C}_0,\{\bar{C}_j,d_j\}_{j=1}^J\}$. Note that $\Delta\in V^c:=\textbf{X}_\Delta\setminus V$ because $\Delta$ is a costless cemetery in this MDP.

\begin{proposition}\label{15JuneProposition02}
Suppose Condition \ref{ConstrainedPPZcondition01} is satisfied.
\begin{itemize}
\item[(a)]
The function $x\in\textbf{X}_\Delta\rightarrow \inf_{\pi}E_x^\pi\left[ \sum_{n=0}^\infty \sum_{j=0}^J \bar{C}_j(X_n,B_{n+1}) \right]$ is lower semicontinuous, so that in particular, the set $V^c$ is a closed subset of $\textbf{X}_\Delta.$
\item[(b)] There is a deterministic stationary strategy, identified by a measurable mapping $f^\ast$ from $\textbf{X}_\Delta$ to $\textbf{B}$, such that
$\inf_{\pi}E_x^\pi\left[ \sum_{n=0}^\infty \sum_{j=0}^J \bar{C}_j(X_n,B_{n+1}) \right]= E_x^{f^\ast}\left[ \sum_{n=0}^\infty  \sum_{j=0}^J \bar{C}_j(X_n,B_{n+1}) \right]$
for each $x\in\textbf{X}_\Delta,$ and
\begin{eqnarray}\label{15JuneEqn08}
Q(V|x,f^\ast(x))=0~\forall~ x\in V^c.
\end{eqnarray}
\end{itemize}
\end{proposition}
\par\noindent\textit{Proof.}  This statement actually concerns an unconstrained version of the impulse control problem, which is the object studied in \cite{PiunovskiySasha:2018}. Both parts (a) and (b) follow from \cite[Thm.1]{PiunovskiySasha:2018}, whereas the last assertion in part (b) is valid further by \cite[Prop.3.2]{Dufour:2012} and its proof. $\hfill\Box$\bigskip

By the definition of the set $V^c$, we see that the deterministic stationary strategy $f^\ast$ coming from Proposition \ref{15JuneProposition02}(b) satisfies
\begin{eqnarray}\label{15JuneEqn07}
\sum_{j=0}^J \bar{C}_j(x,f^\ast(x))=0~\forall~x\in V^c.
\end{eqnarray}

On the set $V^c$, it is intuitively clear that one should apply the strategy $f^\ast$ from Proposition \ref{15JuneProposition02} because it leads to null cost thereafter. On the other hand, on the set $V\times \textbf{B}$, any feasible measure with a finite value in the linear program (\ref{SashaLp01}) possesses some desired properties as stated in the next proposition.

\begin{proposition}\label{15JuneProposition01}
Suppose Conditions \ref{ConstrainedPPZcondition05} and \ref{ConstrainedPPZcondition01} are satisfied.
Then each feasible measure $\mu$ with a finite value in the linear program (\ref{SashaLp01}) is $\sigma$-finite on $V\times \textbf{B}$, and there is a
stochastic kernel $\varphi_\mu$ on $\textbf{B}$ given $V$ satisfying
\begin{eqnarray}\label{15JuneEqn01}
\varphi_\mu(db|x)\mu(dx\times\textbf{B})=\mu(dx\times db)~\mbox{on ${\cal B}(V\times \textbf{B})$.}
\end{eqnarray}
 Moreover,
the stationary strategy $\widetilde{\pi_\mu}$ defined by
\begin{eqnarray}\label{SashaEqn02}
\widetilde{\pi_\mu}(db|x):=\II\{x\in V\}\varphi_\mu(db|x)+\II\{x\in V^c\}\delta_{f^\ast(x)}(db)
\end{eqnarray}
satisfies that
\begin{eqnarray}\label{15JuneEqn06}
\mu^{\widetilde{\pi_\mu}}(\Gamma\times{\bf B})\le \mu(\Gamma\times{\bf B})~\forall~\Gamma\in{\cal B}(V).
\end{eqnarray}
\end{proposition}

\par\noindent\textit{Proof.} The $\sigma$-finiteness of $\mu$ on $V\times \textbf{B}$ can be shown as in the proof of \cite[Thm.3.2]{Dufour:2012}. Then the existence of the stochastic kernel $\varphi_\mu$ is a consequence of this and \cite[Appendix 4]{Dynkin:1979}. The remaining assertions can be shown as in the proof of \cite[Thm.3.3, Cor.3.1]{Dufour:2012}. $\hfill\Box$

\begin{remark}\label{15JuneRemark02}
It is a consequence of (\ref{15JuneEqn07}) and (\ref{15JuneEqn06}) that, in the notation of Proposition \ref{15JuneProposition01}, the occupation measure $\mu^{\widetilde{\pi_\mu}}$ outperforms $\mu$ in the linear program (\ref{SashaLp01}). The verification is as follows:
\begin{eqnarray*}
&&{\cal V}_j(x_0,\pi)=\int_{\textbf{X}_\Delta\times \textbf{B}}\bar{C}_j(x,b)\mu^{\widetilde{\pi_\mu}}(dx\times db)\\
&=&\int_{V\times \textbf{B}}\bar{C}_j(x,b)\widetilde{\pi_\mu}(db|x)\mu^{\widetilde{\pi_\mu}}(dx\times \textbf{B})+\int_{V^c\times \textbf{B}}\bar{C}_j(x,b)\widetilde{\pi_\mu}(db|x)\mu^{\widetilde{\pi_\mu}}(dx\times \textbf{B})\\
&=&\int_{V\times \textbf{B}}\bar{C}_j(x,b)\widetilde{\pi_\mu}(db|x)\mu^{\widetilde{\pi_\mu}}(dx\times \textbf{B})
\le \int_{V\times \textbf{B}}\bar{C}_j(x,b)\widetilde{\pi_\mu}(db|x)\mu (dx\times \textbf{B})\\
&=&\int_{V\times \textbf{B}}\bar{C}_j(x,b)\mu (dx\times db)\le \int_{\textbf{X}_\Delta\times \textbf{B}}\bar{C}_j(x,b)\mu (dx\times db)
\end{eqnarray*}
where the second equality holds by Remark \ref{15JuneRemark01} because $\widetilde{\pi_\mu}$ is a stationary strategy, the third equality follows from (\ref{15JuneEqn07}) and (\ref{SashaEqn02}), the first inequality is by (\ref{15JuneEqn06}), the forth equality is by (\ref{15JuneEqn01}) and (\ref{SashaEqn02}), and the last inequality holds because $\bar{C}_j$ is nonnegative. Consequently, under Conditions \ref{ConstrainedPPZcondition05} and \ref{ConstrainedPPZcondition01},
\begin{eqnarray*}
\mbox{optimal value of problem (\ref{PZZeqn02})}\le ~\mbox{optimal value of the linear program (\ref{SashaLp01})}.
\end{eqnarray*}
\end{remark}

We draw immediately from Remark \ref{15JuneRemark02} and (\ref{15JuneEqn05}) the next corollary.
\begin{corollary}\label{15JuneCorollary01}
Suppose Conditions \ref{ConstrainedPPZcondition05} and \ref{ConstrainedPPZcondition01} are satisfied. The following assertions hold.
 \begin{itemize}
\item[(a)]
The optimal value of problem (\ref{PZZeqn02}) coincides with the optimal value of the linear program (\ref{SashaLp01}).
\item[(b)] The occupation measure of an optimal strategy for problem (\ref{PZZeqn02}) is an optimal solution to the linear program (\ref{SashaLp01}).
    \end{itemize}
\end{corollary}

\begin{remark}\label{rema} 
Every occupation measure as in Definition \ref{15JuneDef02} gives rise to the following `aggregated occupation measure' $\eta$ on ${\bf X}\times{\bf A}_\Box$, where ${\bf A}_\Box:={\bf A}\cup\{\Box\}$ and $\Box$ is an isolated artificial point:
\begin{eqnarray*}
\eta(dy\times\Box)&\defi& \int_{\RR^0_+}\int_{\bf X}\delta_{\phi(x,u)}(dy) \mu(dx\times[u,\infty]\times {\bf A}) du,\\
\eta(dy\times da)&\defi&
 \int_{\bf X}\int_{\RR^0_+} \delta_{\phi(x,\theta)}(dy)  \mu(dx\times d\theta\times da).
\end{eqnarray*}
One can formulate and investigate the linear program, associated with the optimal control problem (\ref{PZZeqn02}), also in terms of aggregated occupation measures. Such linear programs have much in common with works \cite{b13,b13p,b13pp}. This linear programming approach is presented in the follow-up paper \cite{PiunovskiyZhang:2020II}.
\end{remark}

\subsection{A simpler linear program and the linear programming approach}
In view of Proposition \ref{15JuneProposition01}, Remark \ref{15JuneRemark02} and Corollary \ref{15JuneCorollary01}, the main use of the optimal solution say $\mu$ to the linear program (\ref{SashaLp01}) is to produce via (\ref{15JuneEqn01}) the (stationary) strategy to be applied when the state is in $V$. It is natural to ``restrict'' the linear program (\ref{SashaLp01}) to the space of $\sigma$-finite measures on $V\times \textbf{B}$. This results in a simpler linear program: \begin{eqnarray}\label{SashaLp02}
\mbox{Minimize}&:&\int_{V \times\textbf{B}}\bar{C}_0(x,b)\mu(dx\times db) \mbox{ over $\sigma$-finite measures $\mu$ on $V\times \textbf{B}$} \\
 \mbox{suject to}&:& \mu(dx\times \textbf{B})=\delta_{x_0}(dx)+\int_{V\times \textbf{B}}Q(dx|y,b)\mu(dy\times db) \mbox{~on ${\cal B}(V)$};\nonumber\\
&&\int_{V \times\textbf{B}}\bar{C}_j(x,b)\mu(dx\times db)\le d_j,~j=1,2,\dots,J. \nonumber
\end{eqnarray}
(Recall that on $V^c$, the strategy $f^\ast$ is the best one to be used, and it satisfies for all $x\in V^c$ that $\bar{C}_j(x,f^\ast(x))=0$ by the definition of $V^c$, and $Q(V|x,f^\ast(x))=0$  by Proposition \ref{15JuneProposition02}.)

In the next two lemmas, we will elaborate on the relations between the linear programs (\ref{SashaLp01}) and (\ref{SashaLp02}).

\begin{lemma}\label{15JuneLemma01}
Suppose Conditions \ref{ConstrainedPPZcondition05} and \ref{ConstrainedPPZcondition01} are satisfied. Let $\mu'$ be a feasible measure with a finite value in the linear program (\ref{SashaLp01}). Consider the occupation measure $\mu^{\widetilde{\pi_{\mu'}}}$ of the stationary strategy $\widetilde{\pi_{\mu'}}$ defined by (\ref{SashaEqn02}) with $\mu$ being replaced by $\mu'$. Then the restriction of $\mu^{\widetilde{\pi_{\mu'}}}$ on $V\times \textbf{B}$ defined by
\begin{eqnarray*}
\mu^{\widetilde{\pi_{\mu'}}}|_{V\times \textbf{B}}(dx\times db):=\mu^{\widetilde{\pi_{\mu'}}} (dx\times db \cap V\times \textbf{B})
\end{eqnarray*}
is a feasible measure with a finite value in the linear program (\ref{SashaLp02}) satisfying \begin{eqnarray*}
\int_{\textbf{X}_\Delta\times \textbf{B}}\bar{C}_j(x,b)\mu' (dx\times db)\ge \int_{\textbf{X}_\Delta\times \textbf{B}}\bar{C}_j(x,b)\mu^{\widetilde{\pi_{\mu'}}} (dx\times db)=\int_{V\times \textbf{B}}\bar{C}_j(x,b)\mu^{\widetilde{\pi_{\mu'}}}|_{V\times \textbf{B}} (dx\times db).
\end{eqnarray*}
In particular, the optimal value of the linear program (\ref{SashaLp02}) is majorized by the optimal value of the linear program (\ref{SashaLp01}).
\end{lemma}

\par\noindent\textit{Proof.}
For the feasibility of $\mu^{\widetilde{\pi_{\mu'}}}|_{V\times \textbf{B}}$ in the linear program (\ref{SashaLp02}), one may refer to (\ref{15JuneEqn08}), (\ref{15JuneEqn07}) and the definition of $\widetilde{\pi_{\mu'}}(db|x)=\delta_{f^\ast(x)}$ when $x\in V^c$. The rest of this lemma follows from applying Remark \ref{15JuneRemark02} to $\mu'$ and applying Remark \ref{15JuneRemark01} to $\widetilde{\pi_{\mu'}}$.
$\hfill\Box$

\begin{lemma}\label{15JuneLemma02} Suppose Conditions \ref{ConstrainedPPZcondition05} and \ref{ConstrainedPPZcondition01} are satisfied. Let $\mu$ be a feasible measure with a finite value in the linear program (\ref{SashaLp02}), and $\varphi_\mu$ be the stochastic kernel on $\textbf{B}$ given $V$ satisfying (\ref{15JuneEqn01}) and $\widetilde{\pi_\mu}$ be the stationary strategy satisfying (\ref{SashaEqn02}). Then $\mu^{\widetilde{\pi_\mu}}(dx\times{\bf B})\le \mu(dx\times{\bf B})$ on ${\cal B}(V)$, and
\begin{eqnarray*}
\int_{\textbf{X}_\Delta\times \textbf{B}}\bar{C}_j(x,b)\mu^{\widetilde{\pi_\mu}}(dx\times db)=\int_{V\times \textbf{B}}\bar{C}_j(x,b)\mu^{\widetilde{\pi_\mu}}(dx\times db) \le \int_{V\times \textbf{B}}\bar{C}_j(x,b)\mu(dx\times db).
\end{eqnarray*}
In particular, the optimal value of the linear program (\ref{SashaLp02}) is greater or equal to the optimal value of the linear program (\ref{SashaLp01})
\end{lemma}

\par\noindent\textit{Proof.}    The first assertion can be shown as in the proof of \cite[Thm.3.3]{Dufour:2012}. The second assertion follows from the first assertion, just as in Remark \ref{15JuneRemark01}. The last assertion immediately follows from the second assertion. $\hfill\Box$
\bigskip

We now combine Lemmas \ref{15JuneLemma01} and \ref{15JuneLemma02} for the next corollary, which establishes the linear programming approach for solving the impulse control problem (\ref{PZZeqn02}) and will be referred to in the proof of Theorem \ref{Theorem01} in the next section.
\begin{corollary}\label{15JuneCorollary02}
Suppose Conditions \ref{ConstrainedPPZcondition05} and \ref{ConstrainedPPZcondition01} are satisfied. Then the following assertions hold.
\begin{itemize}
\item[(a)] The optimal values of the linear programs (\ref{SashaLp01}) and (\ref{SashaLp02}) as well as the impulse control problem (\ref{PZZeqn02}) are all the same.

\item[(b)] Suppose the measure $\mu_\ast$ is optimal in the linear program (\ref{SashaLp02}). Then   $\mu^{\widetilde{\pi_{\mu_\ast}}}$ defined by (\ref{SashaEqn02}) with $\mu$ being replaced by $\mu_\ast$ is an optimal solution to the linear program (\ref{SashaLp01}), which is feasible with a finite value. In particular, the stationary strategy $\widetilde{\pi_{\mu_\ast}}$ is optimal for the impulse control problem (\ref{PZZeqn02}).
 \item[(c)] If $\mu'$ is an optimal solution to the linear program (\ref{SashaLp01}), then there is an optimal solution to the linear program (\ref{SashaLp02}).
     \end{itemize}
\end{corollary}
\par\noindent\textit{Proof.} By the last assertions of Lemma \ref{15JuneLemma01} and Lemma \ref{15JuneLemma02}, the optimal value of the linear program (\ref{SashaLp02}) is equal to the optimal value of the linear program (\ref{SashaLp01}). Corollary \ref{15JuneCorollary01}(a) asserts that the optimal value of  the linear program (\ref{SashaLp01}) is the same as the optimal value of the impulse control problem (\ref{PZZeqn02}). Part (a) is thus verified. Parts (c) and (b) immediately follow from part (a) and Lemmas \ref{15JuneLemma01} and \ref{15JuneLemma02}, respectively. $\hfill\Box$

\section{Optimality result}\label{sec3}
In this section, we show the existence of a stationary optimal strategy for the impulsive control problem (\ref{PZZeqn02}), as well as the solvability of the linear program (\ref{SashaLp02}). To this end, we hope to apply the relevant results for MDPs. It will be explained that this application is not immediate, but plausible after additional arguments.

\begin{theorem}\label{Theorem01}
Suppose Conditions \ref{ConstrainedPPZcondition05} and \ref{ConstrainedPPZcondition01} are satisfied. Then the following assertions hold.
\begin{itemize}
\item[(a)] There exists an optimal solution say $\mu^\ast$ to the linear program (\ref{SashaLp02}), and an optimal stationary strategy $\widetilde{\pi}$ for the impulse control problem (\ref{PZZeqn02}) such that
    \begin{eqnarray*}
    \mu^\ast(dx\times db)=\mu^\ast(dx\times \textbf{B})\widetilde{\pi}(db|x)~\mbox{on ${\cal B}(V\times \textbf{B})$}
    \end{eqnarray*} and
    $\widetilde{\pi}(db|x)=\delta_{f^\ast(x)}(db)$ for $x\in V^c$, where $f^\ast$ is a deterministic stationary strategy coming from Proposition \ref{15JuneProposition02}
\item[(b)] The linear program (\ref{SashaLp02}) has the same optimal value as the impulse control problem (\ref{PZZeqn02}). Moreover, if $\pi^\ast$ is a strategy such that the restriction of its occupation measure $\mu^{\pi^\ast}|_{V\times \textbf{B}}$ on $V\times \textbf{B}$ solves the linear program (\ref{SashaLp02}), and for each $x_i\in V^c$, $\pi^\ast_{i+1}(db|h_i)$ is concentrated on $\{b\in\textbf{B}:~\sum_{j=0}^J\bar{C}_j(x_i,b)=0\}$, then $\pi^\ast$ is an optimal strategy for the impulse control problem (\ref{PZZeqn02}).
\end{itemize}
\end{theorem}

\par\noindent\textit{Proof.} (a) If we can prove that the linear program (\ref{SashaLp01}) has an optimal solution say $\mu'$, then we may refer to Corollary \ref{15JuneCorollary02}(c) for an optimal solution $\mu_\ast$ to the linear program (\ref{SashaLp02}), and Corollary \ref{15JuneCorollary02}(b) for the statement to be proved with $\mu^\ast=\mu^{\widetilde{\pi_{\mu_\ast}}}$ and $\widetilde{\pi}=\widetilde{\pi_{\mu_\ast}}$. (Recall from Remark \ref{15JuneRemark01} that $\mu^{\widetilde{\pi}}(dx\times db)=\mu^{\widetilde{\pi}}(dx\times\textbf{B})\widetilde{\pi}(db|x)$ on ${\cal B}(\textbf{X}_\Delta\times \textbf{B})$ for any stationary strategy $\widetilde{\pi}$.)

The rest of the proof of part (a) verifies that the linear program (\ref{SashaLp01}) has an optimal solution. To this end, by Corollary \ref{15JuneCorollary01}(b), it suffices to show that the impulse control problem (\ref{PZZeqn02}) has an optimal strategy. Since problem (\ref{PZZeqn02}) can be regarded as an MDP, we would like to do this by making use of a relevant result of MDPs, though it is not immediately applicable as to be seen below.

Let us quote the relevant fact of MDPs from \cite[Thm.4.1]{Dufour:2012}. Namely, for a constrained total cost MDP\footnote{Here the system primitives of the MDP are generic, and are not necessarily defined using the primitives of the impulse control model.} with Borel state space $\textbf{X}_\Delta$, Borel action space $\textbf{B}$, transition probability $Q$, and positive cost functions $\{\bar{C}_j\}_{j=0}^J$, if the model is semicontinuous, then, provided that there exists a feasible strategy with finite value, there is an optimal stationary strategy. Here the model is called semicontinuous if its action space $\textbf{B}$ is compact, $\{\bar{C}_j\}_{j=0}^J$ are all lower semicontinuous, and $Q$ is continuous, i.e., for each bounded continuous function $f$ on $\textbf{X}_\Delta$, $\int_{\textbf{X}_\Delta}f(y)Q(dy|x,a)$ is continuous on $\textbf{X}_\Delta\times \textbf{B}.$

Therefore, under Condition \ref{ConstrainedPPZcondition05}, if the MDP model $\{\textbf{X}_\Delta,\textbf{B},Q,\bar{C}_0,\{\bar{C}_j,d_j\}_{j=1}^J\}$  defined in Section \ref{15JuneSec02} using the primitives of the impulse control model were semicontinuous, we would then be able to refer directly to the quoted result to complete the proof of part (a). It is just that under Condition \ref{ConstrainedPPZcondition05} and \ref{ConstrainedPPZcondition01}, this MDP is not necessarily semicontinuous, however.

In greater detail, under Condition \ref{ConstrainedPPZcondition01}, while the action space $\textbf{B}$ is compact, and the cost functions are lower semicontinuous (as verified in the proof of \cite[Thm.1]{PiunovskiySasha:2018}), the transition probability $Q$ is in general not continuous, because, when $x\ne\Delta$, $\theta_n\in[0,\infty)$ and $\theta_n\to \infty$, the probabilities $Q(dy|x,(\theta_n,a))$ do not converge to $\delta_\Delta(dy)$ in the standard weak topology on the space of Borel probability measures on $\textbf{X}_\Delta$ (generated by bounded continuous functions).

We will get over this difficulty by introducing an auxiliary MDP model $\{\hat{\textbf{X}},\textbf{B},\hat{Q}, \hat{C}_0, \{\hat{C}_j,d_j\}_{j=1}^J\}$. Here the state space in the auxiliary model is \begin{eqnarray*}
\hat{\textbf{X}}:=  ([0,\infty)\times\textbf{X}) \cup \{(\infty,\Delta)\}.
\end{eqnarray*}
We endow it with the metric $\rho_{\hat{\textbf{X}}}$ introduced in Definition \ref{15JuneDef01}, see also Lemma \ref{15JuneLemma05}. By Lemma \ref{15JuneLemma06}, this metric space $\hat{\textbf{X}}$ is a Borel space. The action space $\textbf{B}$ is the same as in the original MDP model.

According to the last assertion of Lemma \ref{15JuneLemma05}, with this metric $\rho_{\hat{\textbf{X}}}$ , $(s_n,x_n)\rightarrow (s,x)\in [0,\infty)\times \textbf{X}$ if and only if $s_n\rightarrow s$ and $x_n\rightarrow x$, whereas if $(s,x)=(\infty,\Delta)$, then the convergence takes place if and only if $s_n\rightarrow \infty$ in the extended half line $[0,\infty]$. Therefore, the relative topology on $[0,\infty)\times \textbf{X}\subseteq \hat{\textbf{X}}$ is the product topology. In that sense we understand the Borel $\sigma$-algebra ${\cal B}([0,\infty)\times \textbf{X}).$ Moreover, by \cite[Lem.7.4]{Bertsekas:1978}, ${\cal B}(\hat{\textbf{X}})$ is generated by ${\cal B}([0,\infty)\times \textbf{X})$ and $\{(\infty,\Delta)\}$. Consequently, we legitimately define the transition kernel $\hat{Q}$ on $\hat{\textbf{X}}$ given $\hat{\textbf{X}}\times\textbf{B}$ by
\begin{eqnarray*}
 \hat{Q}(\Gamma_1\times\Gamma_2|(s,x),(\theta,a)):=\delta_{s+\theta}(\Gamma_1) Q(\Gamma_2|x,(\theta,a))
\end{eqnarray*}
for each $\Gamma_1\in {\cal B}([0,\infty))$ and $\Gamma_2\in{\cal B}(\textbf{X}),$ and \begin{eqnarray*}
\hat{Q}(\{\infty,\Delta\}|(s,x),(\theta,a)):=\II\{\theta+s=\infty\}.
\end{eqnarray*}
Finally, the cost functions are
\begin{eqnarray*}
\hat{C}_j((s,x),(\theta,a)):=\bar{C}_j(x,(\theta,a))
\end{eqnarray*}
for each $(s,x)\in\hat{\textbf{X}}$ and $(\theta,a)\in \textbf{B}$, and $j=0,1,\dots,J.$

Compared to the original MDP model, in this auxiliary MDP model, whose primitives are often signified with a ``hat'', the state space has been extended and a state has an additional coordinate. More precisely, the second coordinate of the current state in the auxiliary model records the current state in the original model, and the first coordinate records the time moment when the impulse was applied, leading to that recorded state in the original model. Therefore,
we take the initial distribution in the hat MDP model as the product measure $\delta_0(dt)\times \delta_{x_0}(dx)$.

Since the consecutive time moments of impulses can be calculated by summing up the time durations between consecutive impulses, i.e., the first coordinate of the action in the original MDP model, any strategy in the auxiliary model admits an equivalent strategy in the original model. Therefore, under Condition \ref{ConstrainedPPZcondition05}, the auxiliary MDP has a feasible strategy with finite value. Furthermore, if there is an optimal strategy in the auxiliary model, then the original MDP admits an optimal strategy, too.

Thus, to complete this proof, it suffices to show that the auxiliary MDP has an optimal strategy. This is done by verifying that under Condition \ref{ConstrainedPPZcondition01} the auxiliary MDP model is semicontinuous, as follows.

The action space $\textbf{B}$ is compact trivially.

For the continuity of $\hat{Q}$, consider a bounded continuous function $f$ on $\hat{\textbf{X}}$. Let $(s_n,x_n)\rightarrow (s,x)$, and $(\theta_n,a_n)\rightarrow (\theta,a).$
Then consider \begin{eqnarray}\label{SashaEqn01}
\int_{\hat{\textbf{X}}}f(t,y)\hat{Q}(dt\times dy|(s_n,x_n),(\theta_n,a_n))&=&f(s_n+\theta_n,l(\phi(x_n,\theta_n),a_n))\II\{s_n+\theta_n< \infty\}\nonumber\\
&&+\II\{s_n+\theta_n= \infty)\}f(\infty,\Delta),
\end{eqnarray}
If $s+\theta<\infty$, then $s_n,\theta_n<\infty$, and $x_n\in\textbf{X}$ for all large enough $n\ge 1,$ and the right hand side of the above equality converges to $f(s+\theta,l(\phi(x,\theta),a))=\int_{\hat{\textbf{X}}}f(t,y)\hat{Q}(dt\times dy|(s,x),(\theta,a))$ according to Condition \ref{ConstrainedPPZcondition01}(b,c). If $s=\infty$ or $\theta=\infty,$ then $s_n+\theta_n\rightarrow \infty$, and hence $f(s_n+\theta_n,l(\phi(x_n,\theta_n),a_n))\rightarrow f(\infty,\Delta)$ according to the definition of the topology on $\hat{\textbf{X}}$. Thus, the right hand side of (\ref{SashaEqn01}) still converges to $f(\infty,\Delta)=\int_{\hat{\textbf{X}}}f(t,y)\hat{Q}(dt\times dy|(s,x),(\theta,a))$, as required.

For the lower semicontinuity of $\hat{C}_j$, where $j=0,1,\dots,J$ is fixed, consider $(s_n,x_n)\rightarrow (s,x)$, and $(\theta_n,a_n)\rightarrow (\theta,a).$ If $(s,x)\in[0,\infty)\times \textbf{X}$, then $(s_n,x_n)\in[0,\infty)\times \textbf{X}$ for all large enough $n\ge 1,$ and so \begin{eqnarray*}
&&\liminf_{n\rightarrow \infty} \hat{C}_j((s_n,x_n),(\theta_n,a_n))=\liminf_{n\rightarrow \infty}\left\{\int_{[0,\theta_n]}  C^g_j(\phi(x_n,u)) du
+ \II\{\theta_n<+\infty\} C^I_j(\phi(x_n,\theta_n),a_n)\right\}\\
&\ge&\int_{[0,\theta]}  C^g_j(\phi(x,u)) du
+ \II\{\theta<+\infty\} C^I_j(\phi(x,\theta),a)= \hat{C}_j((s,x),(\theta,a)),
\end{eqnarray*}
where the inequality holds because the function
$(\theta,x)\in[0,\infty]\times \textbf{X}\rightarrow \int_{[0,\theta]}  C^g_j(\phi(x,u)) du
+ \II\{\theta<\infty\} C^I_j(\phi(x,\theta),a)$ is lower semicontinuous under Condition \ref{ConstrainedPPZcondition01}(c,d,e), as was proved in the proof of \cite[Thm.1]{PiunovskiySasha:2018}, see equation (3.3) therein. If $(s,x)=(\infty,\Delta)$, then $\liminf_{n\rightarrow \infty} \hat{C}_j((s_n,x_n),(\theta_n,a_n))\ge 0=\hat{C}_j((\infty,\Delta),(\theta,a))$ still holds. The lower semicontinuity of $\hat{C}_j$ is now seen. This completes the proof of part (a).

(b) By Corollaries \ref{15JuneCorollary01} and \ref{15JuneCorollary02}, the linear programs (\ref{SashaLp01}) and (\ref{SashaLp02}) have the same optimal value as the impulse control problem (\ref{PZZeqn02}). Now the statement to be proved follows from the observation that $\mu^{\pi^\ast}$ is an optimal solution to the linear program (\ref{SashaLp01}), which holds because $\int_{V^c\times \textbf{B}}\bar{C}_j(x,b)\mu^{\pi^\ast}(dx\times db)=0$.
\hfill $\Box$

\section{Conclusion}\label{15JuneSectionConclusion}
To sum up, under a natural set of compactness-continuity conditions on the system primitives, we  established a linear programming approach to the concerned deterministic optimal impulse control problem, and proved the existence of a stationary optimal control strategy out of a more general class of randomized strategies. This was done by exploiting the relevant facts of MDPs, which were not directly applicable but became plausible after additional arguments.

The linear programming approach proved to be useful for many important control problems, especially in case of constrained optimization. It looks promising to develop such an approach also for non-impulsive deterministic optimal control problems with constraints and to understand its connection with the maximum principle, in case the flow comes from a controlled ordinary differential equation.

\appendix
\section{A metric on the state space of the auxiliary model}
In this appendix, we introduce with justifications a metric $\rho_{\hat{\textbf{X}}}$ on the state space $\hat{\textbf{X}}$ of the auxiliary MDP model in the proof of Theorem \ref{Theorem01}, and verify that it topologizes $\hat{\textbf{X}}$ to a Borel space.

Let $\rho_{\textbf{X}}$ be a compatible metric on $\textbf{X}.$ Since $\textbf{X}$ is a Borel space, according to the lemma of Urysohn, see \cite[Prop.7.2]{Bertsekas:1978}, it is without loss of generality to assume that $\rho_{\textbf{X}}(x,y)\le 2$ for each $x,y\in\textbf{X}.$
\begin{definition}\label{15JuneDef01}
Define a function $\rho_{\hat{\textbf{X}}}:\hat{\textbf{X}}\times \hat{\textbf{X}}\rightarrow [0,\infty)$ as follows. For each $(x_1,s_1),(x_2,s_2)\in\hat{\textbf{X}}$,
\begin{eqnarray*}
\rho_{\hat{\textbf{X}}}((s_1,x_1),(s_2,x_2)):=\rho_{\textbf{X}}(x_1,x_2) (1-g(s_1)\vee g(s_2))+|g(s_1)-g(s_2)|,
\end{eqnarray*}
where
\begin{eqnarray*}
g(s):=\frac{1}{1+\frac{1}{s}}
\end{eqnarray*}
defines a one-to-one correspondence between $[0,\infty]$ and $[0,1]$, accepting $g(0):=0$ and $g(\infty):=1$.
(Note that in the previous definition, if $x_i=\Delta$ for some $i=1,2$, then $\rho_{\textbf{X}}(x_1,x_2)$ is undefined, but $1-g(s_1)\vee g(s_2)=0$; in this case we formally regard $\rho_{\textbf{X}}(x_1,x_2) (1-g(s_1)\vee g(s_2)):=0.$
This convention should be kept in mind below.)
\end{definition}

\begin{lemma}\label{15JuneLemma05}
$\rho_{\hat{\textbf{X}}}$ is a metric on $\hat{\textbf{X}}.$
With this metric, $(s_n,x_n)\rightarrow (s,x)\in [0,\infty)\times \textbf{X}$ if and only if $s_n\rightarrow s$ and $x_n\rightarrow x$, whereas if $(s,x)=(\infty,\Delta)$, then the convergence takes place if and only if $s_n\rightarrow \infty$ in the Euclidean topology. In particular, the relative topology of $\hat{\textbf{X}}$ on $[0,\infty)\times \textbf{X}$ is the same as the product topology on $[0,\infty)\times \textbf{X}$.

\end{lemma}
\par\noindent\textit{Proof.}  We verify that $\rho_{\hat{\textbf{X}}}$ is indeed a metric on $\hat{\textbf{X}}$ as follows. Evidently, it is $[0,\infty)$-valued and symmetric. Moreover, if $(s_1,x_1)=(s_2,x_2)$, then $\rho_{\hat{\textbf{X}}}((s_1,x_1),(s_2,x_2))=0.$
Let us now show that the opposite direction holds.
For notational convenience, below, we write
\begin{eqnarray}\label{15JuneEqn02}
g(s)=\widetilde{s}
\end{eqnarray} for all $s\in[0,\infty].$
Suppose $\rho_{\hat{\textbf{X}}}((s_1,x_1),(s_2,x_2))=0$. Then $\widetilde{s}_1=\widetilde{s}_2$ necessarily. We distinguish two possibilities.
\begin{itemize}
\item If $\widetilde{s}_1=\widetilde{s}_2<1$, then $\rho_{\textbf{X}}(x_1,x_2)=0$, i.e., $(s_1,x_1)=(s_2,x_2)$.
\item  If $\widetilde{s}_1=\widetilde{s}_2=1$, then necessarily $({s}_1,x_1)=({s}_2,x_2)=(\infty,\Delta),$ as required.
\end{itemize}

It remains to verify that $\rho_{\hat{\textbf{X}}}$ satisfies the triangle inequality. Let $(s_i,x_i)\in\hat{\textbf{X}}$, $i=1,2,3,$ be arbitrarily fixed. In case $\widetilde{s}_1\le\widetilde{s}_2\le \widetilde{s}_3,$
\begin{eqnarray*}
&&\rho_{\hat{\textbf{X}}}((s_1,x_1),(s_3,x_3))+\rho_{\hat{\textbf{X}}}((s_3,x_3),(s_2,x_2))-\rho_{\hat{\textbf{X}}}((s_1,x_1),(s_2,x_2))\\
&=&\rho_{\textbf{X}}(x_1,x_3)(1-\widetilde{s}_3)+\widetilde{s}_3-\widetilde{s}_1+\rho_{\textbf{X}}(x_3,x_2)(1-\widetilde{s}_3)+\widetilde{s}_3-\widetilde{s}_2-\rho_{\textbf{X}}(x_1,x_2)(1-\widetilde{s}_2)-\widetilde{s}_2+\widetilde{s}_1\\
&=&\rho_{\textbf{X}}(x_1,x_3)(1-\widetilde{s}_3) +\rho_{\textbf{X}}(x_3,x_2)(1-\widetilde{s}_3)+2(\widetilde{s}_3-\widetilde{s}_2)-\rho_{\textbf{X}}(x_1,x_2)(1-\widetilde{s}_2)\\
&\ge& (1-\widetilde{s}_3)\rho_{\textbf{X}}(x_1,x_2)-\rho_{\textbf{X}}(x_1,x_2)(1-\widetilde{s}_2)+2(\widetilde{s}_3-\widetilde{s}_2)\\
&=&\rho_{\textbf{X}}(x_1,x_2)(\widetilde{s}_2-\widetilde{s}_3)+2(\widetilde{s}_3-\widetilde{s}_2)\ge 0.
\end{eqnarray*}
In case $\widetilde{s}_2\le\widetilde{s}_1\le \widetilde{s}_3,$
\begin{eqnarray*}
&&\rho_{\hat{\textbf{X}}}((s_1,x_1),(s_3,x_3))+\rho_{\hat{\textbf{X}}}((s_3,x_3),(s_2,x_2))-\rho_{\hat{\textbf{X}}}((s_1,x_1),(s_2,x_2))\\
&=&\rho_{\textbf{X}}(x_1,x_3)(1-\widetilde{s}_3)+\widetilde{s}_3-\widetilde{s}_1+\rho_{\textbf{X}}(x_3,x_2)(1-\widetilde{s}_3)+\widetilde{s}_3-\widetilde{s}_2-\rho_{\textbf{X}}(x_1,x_2)(1-\widetilde{s}_1)-\widetilde{s}_1+\widetilde{s}_2\\
&\ge& \rho_{\textbf{X}}(x_1,x_2)(\widetilde{s}_1-\widetilde{s}_3)-2(\widetilde{s}_1-\widetilde{s}_3)\ge 0.
\end{eqnarray*}
Finally, in case $\widetilde{s}_2\le\widetilde{s}_3\le \widetilde{s}_1,$
\begin{eqnarray*}
&&\rho_{\hat{\textbf{X}}}((s_1,x_1),(s_3,x_3))+\rho_{\hat{\textbf{X}}}((s_3,x_3),(s_2,x_2))-\rho_{\hat{\textbf{X}}}((s_1,x_1),(s_2,x_2))\\
&=&\rho_{\textbf{X}}(x_1,x_3)(1-\widetilde{s}_1)+\widetilde{s}_1-\widetilde{s}_3+\rho_{\textbf{X}}(x_3,x_2)(1-\widetilde{s}_3)+\widetilde{s}_3-\widetilde{s}_2-\rho_{\textbf{X}}(x_1,x_2)(1-\widetilde{s}_1)-\widetilde{s}_1+\widetilde{s}_2\\
&=&(\rho_{\textbf{X}}(x_1,x_3)-\rho_{\textbf{X}}(x_1,x_2))(1-\widetilde{s}_1)+\rho_{\textbf{X}}(x_3,x_2)(1-\widetilde{s}_3)\\
&\ge& (\rho_{\textbf{X}}(x_1,x_3)-\rho_{\textbf{X}}(x_1,x_2))(1-\widetilde{s}_1)+\rho_{\textbf{X}}(x_3,x_2)(1-\widetilde{s}_1)\ge  0,
\end{eqnarray*}
as required. The claimed relation holds for all the other possible orders of $\widetilde{s}_1,\widetilde{s}_2,\widetilde{s}_3$ by symmetry. Thus, $\rho_{\hat{\textbf{X}}}$ is a metric on $\hat{\textbf{X}}.$

The last assertion holds automatically from Definition \ref{15JuneDef01}. $\hfill\Box$

\begin{lemma}\label{15JuneLemma06}
Endowed with the metric $\rho_{\hat{\textbf{X}}}$, $\hat{\textbf{X}}$ is a Borel space.
\end{lemma}
\par\noindent\textit{Proof.} We show that $\hat{\textbf{X}}$ is a Borel subset of some Polish (i.e., complete metric) space, as follows. Let $\textbf{X}$ be a Borel subset of a Polish space $[\textbf{X}]$. Consider
\begin{eqnarray*}
\hat{[\textbf{X}]}:=([0,\infty)\times[\textbf{X}]) \cup \{(\infty,\Delta)\},
\end{eqnarray*}
endowed with $\rho_{\hat{[\textbf{X}]}}$ defined as in Definition \ref{15JuneDef01} with $\textbf{X}$ being replaced by $[\textbf{X}]$, see also Lemma \ref{15JuneLemma05} applied to $[\textbf{X}]$.

Let us show that this metric space $\hat{[\textbf{X}]}$ is a Polish space. For its separability, the union of $\{(\infty,\Delta)\}$ and the countable dense subset of $[0,\infty)\times [\textbf{X}]$ (endowed with the product topology) provides a required dense subset of $\hat{[\textbf{X}]},$ by the last assertion of Lemma \ref{15JuneLemma05} applied to $[\hat{\textbf{X}}].$
For its completeness, consider a Cauchy sequence $\{(s_n,x_n)\}_{n=0}^\infty\subseteq \hat{[\textbf{X}]}$. We use the notation (\ref{15JuneEqn02}). Then $\{\widetilde{s}_n\}_{n=0}^\infty$ is Cauchy in $[0,1]$, and so it converges to some $\widetilde{s}\in[0,1]$. If $\widetilde{s}=1,$ then $s_n\rightarrow \infty$ and so $(s_n,x_n)\rightarrow (\infty,\Delta)$ in the metric space $\hat{[\textbf{X}]}$. If $\widetilde{s}\ne 1,$ then $s_n\rightarrow s\ne \infty$, and $x_n\rightarrow x$ for some $x\in [\textbf{X}]$ as $[\textbf{X}]$ is Polish. In either case, Cauchy sequences in $\hat{[\textbf{X}]}$ converge, as required.

Finally, we observe that $\hat{\textbf{X}}$ is a Borel subset of $\hat{[\textbf{X}]}$. In greater detail, $[0,\infty)\times [\textbf{X}]$ is an open subset of the Polish space $\hat{[\textbf{X}]}$, because a singleton in any metric space is closed.
By \cite[Lem.7.4]{Bertsekas:1978}, the trace of ${\cal B}([\hat{\textbf{X}}])$ on $[0,\infty)\times [\textbf{X}]$ is the same as the Borel $\sigma$-algebra ${\cal B}([0,\infty)\times [\textbf{X}])$ when $[0,\infty)\times [\textbf{X}]$ is endowed with the relative topology of $[\hat{\textbf{X}}]$, which is the product topology on $[0,\infty)\times [\textbf{X}]$ according to the last assertion of Lemma \ref{15JuneLemma05} (applied to $[\hat{\textbf{X}}]$). Consequently, $[0,\infty)\times \textbf{X}$ is a Borel subset of $[0,\infty)\times [\textbf{X}]$, and also of $\hat{[\textbf{X}]}$,
and thus $\hat{\textbf{X}}$ is a Borel subset of $\hat{[\textbf{X}]}$. $\hfill\Box$

\section*{Acknowledgements}
This work was partially supported by the Royal Society International Exchanges award IE160503. The authors are grateful to Professor Alexander Plakhov from University of Aveiro (Portugal) and Institute for Information Transmission Problems (Russia) for fruitful discussions.

\end{document}